\documentclass[english,11pt]{article}
\usepackage{amsfonts}
\usepackage[latin1]{inputenc}
\usepackage{amsmath,amsthm,amssymb}
\usepackage{graphicx}
\usepackage{color}

\def\neweq#1{\begin{equation}\label{#1}}
\def\endeq{\end{equation}}

\newtheorem{theorem}{Theorem}[section]

\newtheorem{remark}{Remark}[section]
\input{tcilatex}
\begin{document}

\title{\textbf{Radial solutions for a quasilinear elliptic system of Schrö%
dinger type}}
\author{Dragos-Patru Covei\break \\
{\small \ Constantin Brancusi University of \ Tg-Jiu, Bld. Republicii 1,
210152, Romania}\\
{\small E-mail: c\texttt{ovdra@yahoo.com}}}
\date{}
\maketitle

\begin{abstract}
In this paper we analyze the existence of entire radially symmetric
solutions for Schrodinger system type $\Delta _{p_{i}}u_{i}+h_{i}\left(
r\right) \left\vert \nabla u_{i}\right\vert ^{p_{i}-1}=a_{i}\left( r\right)
f_{i}\left( u_{1},...,u_{d}\right) $ for $i=1,...,d$ on $R^{N}$ where $%
p_{i}>1$, $d\in \{1,2,3,...\}$, $h_{i}$ and $a_{i}$ are nonnegative radial
continuous functions and $f_{i}$ are nonnegative increasing continuous
functions on $\left[ 0,\infty \right) $.\\[3pt]
\end{abstract}

\baselineskip16pt \renewcommand{\theequation}{\arabic{section}.%
\arabic{equation}} \catcode`@=11 \@addtoreset{equation}{section} \catcode%
`@=12



\section{Introduction}

In this article we continue to study the existence results for systems such
as 
\begin{equation}
\left\{ 
\begin{array}{l}
\Delta _{p_{1}}u_{1}+h_{1}\left( r\right) \left\vert \nabla u_{1}\right\vert
^{p_{1}-1}=a_{1}\left( r\right) f_{1}\left( u_{1},...,u_{d}\right) \text{ in 
}\mathbb{R}^{N}\text{ }, \\ 
... \\ 
\Delta _{p_{d}}u_{d}+h_{d}\left( r\right) \left\vert \nabla u_{d}\right\vert
^{p_{d}-1}=a_{d}\left( r\right) f_{d}\left( u_{1},...,u_{d}\right) \text{ in 
}\mathbb{R}^{N},%
\end{array}%
\right.  \label{11}
\end{equation}%
where $r:=\left\vert x\right\vert \geq 0$ denotes the Euclidean length of $%
x\in \mathbb{R}^{N}$ , $N\geq 3$, $d\geq 1$ is integer, $\Delta _{p_{j}}$ ($%
j=1,...,d$) is the $p_{j}$-Laplacian operator defined by 
\begin{equation*}
\Delta _{p_{j}}u:=\func{div}(\left\vert \nabla u\right\vert ^{p_{j}-2}\nabla
u)\text{, }1<p_{j}<\infty ,\text{ }
\end{equation*}
$h_{j}$, $a_{j}:\left[ 0,\infty \right) \rightarrow \left[ 0,\infty \right) $
are radial continuous functions and $f_{j}$ satisfy the following hypotheses

(C1)\quad $f_{j}:\left[ 0,\infty \right) ^{d}\rightarrow \left[ 0,\infty
\right) $ are continuous in all variables;

(C2)\quad $f_{j}$ are non-decreasing on $\left[ 0,\infty \right) ^{d}$ in
all variables.

Particular forms of the Schrodinger system type like (\ref{11}) have been
considered in \cite{GB}-\cite{YE} and references therein. The problem
considered in our paper it is considered in a general form that includes
various types of problems from applied sciences. For example, the
time-independent Schrodinger equation \textit{\ }%
\begin{equation}
(h^{2}/2m)\Delta u=(V-E)u  \label{s}
\end{equation}%
where $h$ is the Plank constant, $m$ is the mass of a particle moving under
the action of a force field described by the potential $V$ whose wave
function is $u$ and the quantity $E$ is the total energy of the particle
(see the book \cite{GG}). Also, some particular classes of systems studied
in our work are used in the description of several physical phenomena such
as the propagation of pulses in birefringent optical fibers and Kerr-like
photorefractive media (see the articles \cite{AA}, \cite{M}). Moreover, in
the absence of nonlinear gradient term such problems appear in medical
science (see \cite{LL}).

Several theoretical results are available in the literature for \ the
problem of the form (\ref{11}). Most of the studies are about the existence
or the nonexistence of positive radial ones, because the applications have
been concentrated in symmetry theory. In particular, our research is closely
related to \cite{CD3,LA2,T1} where the authors have obtained some
theoretical interesting results and the paper \cite{LAS} where the
applications can be found. \ In \cite{T1}, the author consider the \textit{%
Schrödinger} system type%
\begin{equation}
\left\{ 
\begin{array}{l}
\Delta _{p_{1}}u_{1}=H_{1}\left( \left\vert x\right\vert \right)
u_{2}^{\alpha _{1}}\text{ }, \\ 
\Delta _{p_{2}}u_{2}=H_{2}\left( \left\vert x\right\vert \right)
u_{3}^{\alpha _{2}}\text{ }, \\ 
... \\ 
\Delta _{p_{m}}u_{m}=H_{m}\left( \left\vert x\right\vert \right)
u_{m+1}^{\alpha _{m}}\text{ },\text{ }u_{m+1}=u_{1},%
\end{array}%
\right. \text{ }x\in \mathbb{R}^{N},  \label{T}
\end{equation}%
where $\alpha _{i}$ and $p_{i}$ ($i=1,2,...,m$) are constants satisfying 
\begin{equation*}
\alpha _{1}\cdot \cdot \cdot \alpha _{m}>(p_{1}-1)\cdot \cdot \cdot
(p_{m}-1),
\end{equation*}%
$\Delta _{p_{i}}$ ($i=1,2,...,m$) is the $p_{i}$-Laplacian operator and the
functions $H_{i}$ ($i=1,2,...,m$), are nonnegative continuous functions on $%
\left[ 0,\infty \right) $. Under these hypotheses and some additional
conditions he established the existence and non-existence results for
solutions of the system (\ref{T}) which will be improved here. Our solving
method gives a stronger meaning to the obtained solutions, compared with the
famous book \textquotedblleft Particle Physics and the Schrodinger
Equation\textquotedblright\ \cite{GG} where similar solutions are detected.
\ 

Before we start to describe our results let us mention that the problem of
existence of solutions to (\ref{11}) has received an increased interest with
Alan Lair's recent paper \cite{LA2}. In \cite{LA2} the author considered
problem%
\begin{equation}
\left\{ 
\begin{array}{l}
\Delta u_{1}=a_{1}\left( \left\vert x\right\vert \right) u_{2}^{\alpha
}\left( \left\vert x\right\vert \right) \text{ },\text{ }(\alpha \in \left(
0,1\right] ), \\ 
\Delta u_{2}=a_{2}\left( \left\vert x\right\vert \right) u_{1}^{\beta
}\left( \left\vert x\right\vert \right) ,\text{ }(\beta \in \left( 0,1\right]
)\text{ },%
\end{array}%
\right. \text{ }x\in \mathbb{R}^{N}.  \label{l}
\end{equation}%
The main result achieved by Alan Lair may be summarized as follows:

\begin{proposition}
Problem (\ref{l}) has an explosive radial symmetric solution on $\mathbb{R}%
^{N}$ if and only if \textit{continuous radially }symmetric\textit{\
functions } $a_{i_{i=\overline{1,2}}}:\left[ 0,\infty \right) \rightarrow %
\left[ 0,\infty \right) $ simultaneously meet the following conditions 
\begin{eqnarray*}
\int_{0}^{\infty }ta_{1}\left( t\right) \left(
t^{2-N}\int_{0}^{t}s^{N-3}\int_{0}^{s}\tau a_{2}\left( \tau \right) d\tau
ds\right) ^{\alpha }dt &=&\infty \text{ } \\
\int_{0}^{\infty }ta_{2}\left( t\right) \left(
t^{2-N}\int_{0}^{t}s^{N-3}\int_{0}^{s}\tau a_{1}\left( \tau \right) d\tau
ds\right) ^{\beta }dt &=&\infty \text{.}
\end{eqnarray*}
\end{proposition}

Moreover, the author issues the following problem:

"\textit{It remains unknown whether an analogous result for the system}%
\begin{equation}
\left\{ 
\begin{array}{l}
\Delta u_{1}\left( \left\vert x\right\vert \right) =a_{1}\left( \left\vert
x\right\vert \right) f_{1}\left( u_{2}\left( \left\vert x\right\vert \right)
\right) \text{ for }x\in \mathbb{R}^{N}\text{ },\text{ } \\ 
\Delta u_{2}\left( \left\vert x\right\vert \right) =a_{2}\left( \left\vert
x\right\vert \right) f_{2}\left( u_{1}\left( \left\vert x\right\vert \right)
\right) \text{ for }x\in \mathbb{R}^{N}\text{ },%
\end{array}%
\right.  \label{prop}
\end{equation}%
\textit{\ where }$f_{1}$\textit{\ and }$f_{2}$\textit{\ meet, for example
the Keller-Osserman \cite{K,O} condition}%
\begin{equation}
\int_{1}^{\infty }[\int_{0}^{s}f_{i}\left( t\right) dt]^{-1/2}ds=\infty 
\text{, }i=\overline{1,2}\text{,}  \label{r1}
\end{equation}%
\textit{or the Ye and Zhou \cite{YE} condition}%
\begin{equation}
\int_{1}^{\infty }[f_{i}\left( t\right) ]^{-1}dt=\infty \text{, }i=\overline{%
1,2}\text{."}  \label{r2}
\end{equation}

We mention that in the paper \cite{CD2} the author established for the
system (\ref{prop}) a necessary condition as well as a sufficient condition
for a positive radial solution to be large under conditions of the (\ref{r1}%
) type.

Now let us give a detailed description of our results. Throughout this paper
we use the notations 
\begin{equation*}
\begin{array}{ll}
j=1,...,d, & H_{j}\left( r\right) :=r^{N-1}e^{\int_{0}^{r}h_{j}\left(
t\right) dt}, \\ 
A_{j}\left( \infty \right) :=\lim_{r\rightarrow \infty }A_{j}\left( r\right)
, & A_{j}\left( r\right) :=\int_{0}^{r}\left( \frac{1}{H_{j}\left( t\right) }%
\int_{0}^{t}H_{j}\left( s\right) a_{j}\left( s\right) ds\right) ^{\frac{1}{%
p_{j}-1}}dt,%
\end{array}%
\end{equation*}%
and the Lair's (\cite[p.211]{LA3}) quantity%
\begin{equation*}
F\left( \infty \right) =\lim_{r\rightarrow \infty }F\left( r\right) ,\text{ }%
F\left( r\right) =\int_{a}^{r}\left( 1+\underset{j=1}{\overset{d}{\sum }}%
f_{j}\left( s\right) \right) ^{\frac{1}{1-\min \left\{
p_{1},...,p_{d}\right\} }}ds\text{; }r\geq a>0.
\end{equation*}%
We see that 
\begin{equation*}
F^{\prime }\left( r\right) =\left( 1+\underset{j=1}{\overset{d}{\sum }}%
f_{j}\left( r\right) \right) ^{-\frac{1}{\min \left\{
p_{1},...,p_{d}\right\} -1}}>0\text{ for all }r>a
\end{equation*}%
and $F$ has the inverse function $F^{-1}$ on $\left[ a,\infty \right) $.

We now give our main theorems.

\begin{theorem}
\label{1}Suppose that (C1)-(C2) hold and that

(C3)\quad $F\left( \infty \right) =\infty $.

Then the system (\ref{11}) possesses at least one positive radial solution $%
\left( u_{1},...,u_{d}\right) $. If, in addition, $A_{j}\left( \infty
\right) <\infty $ ($j=1,..,d$), the positive radial solution $\left(
u_{1},...,u_{d}\right) $ is bounded. On the other hand, if $A_{j}\left(
\infty \right) =\infty $ the positive solution $\left(
u_{1},...,u_{d}\right) $ is entire large solution, i.e. 
\begin{equation*}
\lim_{r\rightarrow \infty }u_{1}\left( r\right) =...=\lim_{r\rightarrow
\infty }u_{d}\left( r\right) =\infty .
\end{equation*}
\end{theorem}

\begin{theorem}
\label{2} Assume that (C1)-(C2) hold and that

(C4)\quad $F\left( \infty \right) <\infty ;$

(C5)\quad $A_{j}\left( \infty \right) <\infty $ ($j=1,..,d$)$;$

(C6)\quad there exists $\beta >\frac{a}{d}$ such that 
\begin{equation*}
\underset{j=1}{\overset{d}{\sum }}A_{j}\left( \infty \right) <F\left( \infty
\right) -F\left( d\beta \right) .
\end{equation*}%
Then, the system (\ref{11}) possesses at least one positive bounded radial
solution $\left( u_{1},...,u_{d}\right) $ satisfying%
\begin{equation*}
\beta +f_{j}^{1/\left( p_{j}-1\right) }\left( \beta ,...,\beta \right)
A_{j}\left( r\right) \leq u_{j}\left( r\right) \leq F^{-1}\left( F\left(
d\beta \right) +\underset{j=1}{\overset{d}{\sum }}A_{j}\left( r\right)
\right) .
\end{equation*}
\end{theorem}

\begin{theorem}
\textbf{\label{3}} (i) Assume that $A_{i_{i=\overline{1,d}}}\left( \infty
\right) =\infty $ and%
\begin{equation}
\lim_{s\rightarrow \infty }\frac{\underset{i=1}{\overset{d}{\sum }}\left(
1+f_{i}\left( s,...,s\right) \right) ^{\frac{1}{\min \left\{
p_{1},...,p_{d}\right\} -1}}\text{ }}{s}=0.  \label{10}
\end{equation}%
Then the system (\ref{11}) has infinitely many positive entire large
solutions.

(ii) Furthermore, if $A_{i_{i=\overline{1,d}}}\left( \infty \right) <\infty $
and 
\begin{equation*}
\sup_{s\geq 0}\left[ \underset{i=1}{\overset{d}{\sum }}\left( 1+f_{i}\left(
s,...,s\right) \right) ^{\frac{1}{\min \left\{ p_{1},...,p_{d}\right\} -1}}%
\right] <\infty 
\end{equation*}%
then the system (\ref{11}) has infinitely many positive entire bounded
solutions.
\end{theorem}

\section{Proof of Theorems}

\subsection{Proof of Theorem \protect\ref{1}}

We note that radial solutions of system (\ref{11}) are solutions $u_{j}$ ($%
j=1,...,d$) to the ordinary differential system%
\begin{equation}
\left\{ 
\begin{array}{c}
\frac{1}{r^{N-1}}\left( r^{N-1}\left\vert u^{\prime }\right\vert
^{p_{1}-2}u^{\prime }\right) ^{\prime }+h_{1}\left( r\right) \left\vert
u_{1}^{\prime }\right\vert ^{p_{1}-1}=a_{1}\left( r\right) f_{1}\left(
u_{1},...,u_{d}\right) , \\ 
... \\ 
\frac{1}{r^{N-1}}\left( r^{N-1}\left\vert u^{\prime }\right\vert
^{p_{d}-2}u^{\prime }\right) ^{\prime }+h_{d}\left( r\right) \left\vert
u_{d}^{\prime }\right\vert ^{p_{d}-1}=a_{d}\left( r\right) f_{d}\left(
u_{1},...,u_{d}\right) ,%
\end{array}%
\right.  \label{r}
\end{equation}%
and that any solution $\left( u_{1},...,u_{d}\right) $ to the integral
equations 
\begin{equation}
\left\{ 
\begin{array}{l}
u_{1}\left( r\right) =\beta +\int_{0}^{r}\left( \frac{1}{H_{1}\left(
t\right) }\int_{0}^{t}H_{1}\left( s\right) a_{1}\left( s\right) f_{1}\left(
u_{1}\left( s\right) ,...,u_{d}\left( s\right) \right) ds\right) ^{\frac{1}{%
p_{1}-1}}dt,\text{ } \\ 
... \\ 
u_{d}\left( r\right) =\beta +\int_{0}^{r}\left( \frac{1}{H_{d}\left(
t\right) }\int_{0}^{t}H_{d}\left( s\right) a_{d}\left( s\right) f_{d}\left(
u_{1}\left( s\right) ,...,u_{d}\left( s\right) \right) ds\right) ^{\frac{1}{%
p_{d}-1}}dt,%
\end{array}%
\right.  \label{in}
\end{equation}%
is a solution to (\ref{11}).

We will begin by establishing a solution of (\ref{in}) in $\left(
C[0,R]\right) ^{d}$ for arbitrary $R>0$. For this, we apply a standard
iteration procedure by letting 
\begin{equation*}
u_{1}^{0}=...=u_{d}^{0}=\beta >0
\end{equation*}%
the central values for the integral equations system and generating a
non-decreasing sequence $\left\{ u_{j}^{k}\right\} _{1\leq j\leq d}^{k\geq
1} $ in which $u_{j}^{k}$ is calculated from $u_{j}^{k-1}$ by

\begin{equation}
\left\{ 
\begin{array}{l}
u_{1}^{k}\left( r\right) =\beta +\int_{0}^{r}\left( \frac{1}{H_{1}\left(
t\right) }\int_{0}^{t}H_{1}\left( s\right) a_{1}\left( s\right) f_{1}\left(
u_{1}^{k-1},...,u_{d}^{k-1}\right) ds\right) ^{\frac{1}{p_{1}-1}}dt, \\ 
... \\ 
u_{d}^{k}\left( r\right) =\beta +\int_{0}^{r}\left( \frac{1}{H_{d}\left(
t\right) }\int_{0}^{t}H_{d}\left( s\right) a_{d}\left( s\right) f_{d}\left(
u_{1}^{k-1},...,u_{d}^{k-1}\right) ds\right) ^{\frac{1}{p_{d}-1}}dt.%
\end{array}%
\right.  \label{ints}
\end{equation}%
Due to the form (\ref{ints}), for all $r\geq 0,$ $j=\overline{1,d}$ and $%
k\in N$ we have $u_{j}^{k}\left( r\right) \geq \beta $. Furthermore, we can
easy see that $\left\{ u_{j}^{k}\right\} _{1\leq j\leq d}^{k\geq 1}$ are
non-decreasing sequence on $\left[ 0,\infty \right) $.

By conditions (C1) and (C2) we obtain%
\begin{equation}
\begin{array}{l}
\left( u_{1}^{k}\left( r\right) \right) ^{\prime }=\left( \frac{1}{%
H_{1}\left( r\right) }\int_{0}^{r}H_{1}\left( s\right) a_{1}\left( s\right)
f_{1}\left( u_{1}^{k-1}\left( s\right) ,...,u_{d}^{k-1}\left( s\right)
\right) ds\right) ^{\frac{1}{p_{1}-1}} \\ 
\leq f_{1}^{\frac{1}{p_{1}-1}}\left( u_{1}^{k}\left( r\right)
,...,u_{d}^{k}\left( r\right) \right) A_{1}^{\prime }\left( r\right) \leq
\left( \underset{j=1}{\overset{d}{\sum }}f_{j}\left( \overset{d}{\underset{%
j=1}{\sum }}u_{j}^{k}\left( r\right) \right) \right) ^{\frac{1}{p_{1}-1}%
}A_{1}^{\prime }\left( r\right) \\ 
\leq \left( 1+\underset{j=1}{\overset{d}{\sum }}f_{j}\left( \overset{d}{%
\underset{j=1}{\sum }}u_{j}^{k}\left( r\right) \right) \right) ^{\frac{1}{%
\min \left\{ p_{1},...,p_{d}\right\} -1}}A_{1}^{\prime }\left( r\right) , \\ 
... \\ 
\left( u_{d}^{k}\left( r\right) \right) ^{\prime }=\left( \frac{1}{%
H_{d}\left( r\right) }\int_{0}^{r}H_{d}\left( s\right) a_{d}\left( s\right)
f_{d}\left( u_{1}^{k-1}\left( s\right) ,...,u_{d}^{k-1}\left( s\right)
\right) ds\right) ^{\frac{1}{p_{d}-1}} \\ 
\leq f_{d}^{\frac{1}{p_{d}-1}}\left( u_{1}^{k}\left( r\right)
,...,u_{d}^{k}\left( r\right) \right) A_{d}^{\prime }\left( r\right) \leq
\left( \underset{j=1}{\overset{d}{\sum }}f_{j}\left( \overset{d}{\underset{%
j=1}{\sum }}u_{j}^{k}\left( r\right) \right) \right) ^{\frac{1}{p_{d}-1}%
}A_{d}^{\prime }\left( r\right) \\ 
\leq \left( 1+\underset{j=1}{\overset{d}{\sum }}f_{j}\left( \overset{d}{%
\underset{j=1}{\sum }}u_{j}^{k}\left( r\right) \right) \right) ^{\frac{1}{%
\min \left\{ p_{1},...,p_{d}\right\} -1}}A_{d}^{\prime }\left( r\right) .%
\end{array}
\label{sis}
\end{equation}%
Summing up gives%
\begin{equation*}
\left( 1+\underset{j=1}{\overset{d}{\sum }}f_{j}\left( \overset{d}{\underset{%
j=1}{\sum }}u_{j}^{k}\left( t\right) \right) \right) ^{-\frac{1}{\min
\left\{ p_{1},...,p_{d}\right\} -1}}\cdot \left( \overset{d}{\underset{j=1}{%
\sum }}u_{j}^{k}\left( t\right) \right) ^{\prime }\leq \overset{d}{\underset{%
j=1}{\sum }}A_{j}^{\prime }\left( t\right) .
\end{equation*}%
Integrating this over $\left[ 0,r\right] $, produces%
\begin{equation*}
\int_{0}^{r}\left( 1+\underset{j=1}{\overset{d}{\sum }}f_{j}\left( \overset{d%
}{\underset{j=1}{\sum }}u_{j}^{k}\left( t\right) \right) \right) ^{-\frac{1}{%
\min \left\{ p_{1},...,p_{d}\right\} -1}}\cdot \left( \overset{d}{\underset{%
j=1}{\sum }}u_{j}^{k}\left( t\right) \right) ^{\prime }dt\leq \overset{d}{%
\underset{j=1}{\sum }}A_{j}\left( r\right)
\end{equation*}%
for each $r>0$, which can be rewritten as%
\begin{equation*}
\int_{0}^{r}F^{\prime }\left( \overset{d}{\underset{j=1}{\sum }}%
u_{j}^{k}\left( t\right) \right) dt\leq \overset{d}{\underset{j=1}{\sum }}%
A_{j}\left( r\right) \text{ for each }r>0\text{,}
\end{equation*}%
from which we get%
\begin{equation}
F\left( \overset{d}{\underset{j=1}{\sum }}u_{j}^{k}\left( r\right) \right)
-F\left( d\beta \right) \leq \overset{d}{\underset{j=1}{\sum }}A_{j}\left(
r\right) \text{ for all }r\geq 0.  \label{21}
\end{equation}%
Since $F^{-1}$ is increasing on $\left[ 0,\infty \right) $, follows that 
\begin{equation}
\overset{d}{\underset{j=1}{\sum }}u_{j}^{k}\left( r\right) \leq F^{-1}\left(
F\left( d\beta \right) +\overset{d}{\underset{j=1}{\sum }}A_{j}\left(
r\right) \right) \text{ for all }r\geq 0.  \label{22}
\end{equation}%
Since (C3) holds, we can see that%
\begin{equation}
F^{-1}\left( \infty \right) =\infty .  \label{23}
\end{equation}%
It follows that the sequences $\left\{ u_{j}^{k}\right\} _{1\leq j\leq
d}^{k\geq 1}$ \ are bounded and non-decreasing on $\left[ 0,R\right] $ for $%
R>0$.

Thus 
\begin{equation}
\left( u_{1}^{k},...,u_{d}^{k}\right) \text{ converges to }\left(
u_{1},...,u_{d}\right) \text{ on }\left[ 0,R\right] ^{d}.  \label{sub}
\end{equation}%
Consequently $\left( u_{1},...,u_{d}\right) $ is the positive entire radial
solution of system (\ref{11}) in $\overline{B(0,R)}\subset \mathbb{R}^{N}$
with central values $u_{1}\left( 0\right) =...=u_{d}\left( 0\right) =\beta $.

Since $R$ is arbitrary, we can use the diagonal argument to show (\ref{sub})
has a convergent subsequence on $\left( C[0,\infty )\right) ^{d}$ to a
function denoted again by $\left( u_{1},...,u_{d}\right) $ and this is a
solution of (\ref{11}) in $\mathbb{R}^{N}$. A complete proof of this
procedure can be found in the work of (\cite{AF}) where numerical results
are also commented.

In addition, when 
\begin{equation*}
A_{j}\left( \infty \right) <\infty ,\text{ }j=1,...,d
\end{equation*}%
we see by (\ref{22}) that%
\begin{equation*}
\overset{d}{\sum_{j=1}}u_{j}\left( r\right) \leq F^{-1}\left( F\left( d\beta
\right) +\overset{d}{\sum_{j=1}}A_{j}\left( \infty \right) \right) \text{
for all }r\geq 0
\end{equation*}%
when 
\begin{equation*}
A_{j}\left( \infty \right) =\infty \text{ for }j=j=1,...,d
\end{equation*}%
by (C2) and the monotonicity of $\left\{ u_{j}^{k}\right\} _{1\leq j\leq
d}^{k\geq 1}$ follows%
\begin{equation*}
u_{j}\left( r\right) \geq \beta +f_{j}^{1/\left( p_{j}-1\right) }\left(
\beta ,...,\beta \right) A_{j}\left( r\right) \text{, for all }r\geq 0\text{
and }j=1,...,d\text{.}
\end{equation*}%
Then 
\begin{equation*}
\lim_{r\rightarrow \infty }u_{1}\left( r\right) =...=\lim_{r\rightarrow
\infty }u_{d}\left( r\right) =\infty ,
\end{equation*}%
and our proof is complete.

\subsection{Proof of Theorem \protect\ref{2}.}

In a manner similar to our Theorem \ref{1} proof above, we obtain that 
\begin{equation}
F\left( \overset{d}{\underset{j=1}{\sum }}u_{j}^{k}\left( r\right) \right)
\leq F\left( d\beta \right) +\overset{d}{\underset{j=1}{\sum }}A_{j}\left(
\infty \right) <F\left( \infty \right) <\infty .  \label{24}
\end{equation}%
Because $F^{-1}$ is strictly increasing on $\left[ 0,\infty \right) $ we have%
\begin{equation}
\overset{d}{\underset{j=1}{\sum }}u_{j}^{k}\left( r\right) \leq F^{-1}\left(
F\left( d\beta \right) +\overset{d}{\underset{j=1}{\sum }}A_{j}\left( \infty
\right) \right) <\infty \text{ for all }r\geq 0.  \label{25}
\end{equation}%
Moreover, since the sequence $\left\{ u_{j}^{k}\left( r\right) \right\} $ is
monotone it converges to some function $\left\{ u_{j}\left( r\right)
\right\} _{1\leq j\leq d}$ on $\mathbb{R}^{N}$ that in fact is a solution to
(\ref{11}) and the proof is complete.

\subsection{Proof of Theorem \protect\ref{3}}

\bigskip We first see that radial solutions of (\ref{11}) are solutions $%
\left( u_{1},...,u_{d}\right) $ of the differential equations system%
\begin{equation}
\left\{ 
\begin{array}{c}
\left( p_{1}-1\right) \left\vert u_{1}^{\prime }\left( r\right) \right\vert
^{p_{1}-2}u_{1}^{\prime \prime }+\frac{N-1}{r}u_{1}^{\prime }\left( r\right)
^{p_{1}-1}+h_{1}\left( r\right) \left\vert u_{1}^{^{\prime }}\left( r\right)
\right\vert ^{p_{1}-1}=a_{1}\left( r\right) f_{1}\left( u_{1}\left( r\right)
,...,u_{d}\left( r\right) \right) , \\ 
... \\ 
\left( p_{d}-1\right) \left\vert u_{d}^{\prime }\left( r\right) \right\vert
^{p_{d}-2}u_{d}^{\prime \prime }+\frac{N-1}{r}u_{d}^{\prime }\left( r\right)
^{p_{d}-1}+h_{d}\left( r\right) \left\vert u_{d}^{^{\prime }}\left( r\right)
\right\vert ^{p_{d}-1}=a_{d}\left( r\right) f_{d}\left( u_{1}\left( r\right)
,...,u_{d}\left( r\right) \right) .%
\end{array}%
\right.   \label{1010}
\end{equation}%
Since the radial solutions of (\ref{11}) are solutions of the differential
equations system (\ref{1010}) it follows that the radial solutions of (\ref%
{11}) with 
\begin{equation*}
u_{1}(0)=\beta _{1},...,u_{d}(0)=\beta _{d}
\end{equation*}%
where $\beta _{i}$ ($i=1,...,d$)\ may be any non-negative numbers, satisfy: 
\begin{equation}
\left\{ 
\begin{array}{l}
u_{1}\left( r\right) =\beta _{1}+\int_{0}^{r}\left( \frac{1}{H_{1}\left(
t\right) }\int_{0}^{t}H_{1}\left( s\right) a_{1}\left( s\right) f_{1}\left(
u_{1}\left( s\right) ,...,u_{d}\left( s\right) \right) ds\right) ^{1/\left(
p_{1}-1\right) }dt,\text{ } \\ 
... \\ 
u_{d}\left( r\right) =\beta _{d}+\int_{0}^{r}\left( \frac{1}{H_{d}\left(
t\right) }\int_{0}^{t}H_{d}\left( s\right) a_{d}\left( s\right) f_{d}\left(
u_{1}\left( s\right) ,...,u_{d}\left( s\right) \right) ds\right) ^{1/\left(
p_{d}-1\right) }dt.%
\end{array}%
\right.   \label{222}
\end{equation}%
Define 
\begin{equation*}
u_{1}^{0}=\beta _{1},...,u_{d}^{0}=\beta _{d}\text{ for }r\geq 0.
\end{equation*}%
Let $\left\{ u_{j}^{k}\right\} _{j=1,...,d}^{k\geq 1}$ be a sequence of
functions on $\left[ 0,\infty \right) $ given by%
\begin{equation}
\left\{ 
\begin{array}{l}
u_{1}^{k+1}\left( r\right) =\beta _{1}+\int_{0}^{r}\left( \frac{1}{%
H_{1}\left( t\right) }\int_{0}^{t}H_{1}\left( s\right) a_{1}\left( s\right)
f_{1}\left( u_{1}^{k}\left( s\right) ,...,u_{d}^{k}\left( s\right) \right)
ds\right) ^{1/\left( p_{1}-1\right) }dt, \\ 
... \\ 
u_{d}^{k+1}\left( r\right) =\beta _{d}+\int_{0}^{r}\left( \frac{1}{%
H_{d}\left( t\right) }\int_{0}^{t}H_{d}\left( s\right) a_{d}\left( s\right)
f_{d}\left( u_{1}^{k}\left( s\right) ,...,u_{d}^{k}\left( s\right) \right)
ds\right) ^{1/\left( p_{d}-1\right) }dt.%
\end{array}%
\right.   \label{233}
\end{equation}%
We remark that, for all $r\geq 0,$ $j=1,...,d$ and $k\in N$ 
\begin{equation*}
u_{j}^{k}\left( r\right) \geq \beta _{j}\text{.}
\end{equation*}%
Moreover $\left\{ u_{j}^{k}\right\} _{j=1,...,d}^{k\geq 1}$ are
non-decreasing sequence on $\left[ 0,\infty \right) $ such that%
\begin{equation}
u_{i}^{k}\left( r\right) \leq u_{i}^{k+1}\left( r\right) \leq \left[
f_{i}\left( \overset{d}{\underset{j=1}{\Sigma }}u_{j}^{k}\left( r\right)
,...,\overset{d}{\underset{j=1}{\Sigma }}u_{j}^{k}\left( r\right) \right) %
\right] ^{1/\left( p_{i}-1\right) }A_{i}\left( r\right) ,\text{ }i=1,...,d.
\label{26}
\end{equation}%
Let $R>0$ be arbitrary. It is easy to see that (\ref{26}) implies%
\begin{equation*}
\overset{d}{\underset{i=1}{\Sigma }}u_{i}^{k}\left( R\right) \leq \overset{d}%
{\underset{i=1}{\Sigma }}\beta _{i}+\overset{d}{\underset{i=1}{\Sigma }}%
\left[ 1+f_{i}\left( \overset{d}{\underset{j=1}{\Sigma }}u_{j}^{k}\left(
R\right) ,...,\overset{d}{\underset{j=1}{\Sigma }}u_{j}^{k}\left( R\right)
\right) \right] ^{\frac{1}{\min \left\{ p_{1},...,p_{d}\right\} -1}}\overset{%
d}{\underset{i=1}{\Sigma }}A_{i}\left( R\right) ,k\geq 1,
\end{equation*}%
and so%
\begin{equation}
1\leq \frac{\overset{d}{\underset{i=1}{\Sigma }}\beta _{i}}{\overset{d}{%
\underset{i=1}{\Sigma }}u_{i}^{k}\left( R\right) }+\frac{\overset{d}{%
\underset{i=1}{\Sigma }}\left[ 1+f_{i}\left( \overset{d}{\underset{j=1}{%
\Sigma }}u_{j}^{k}\left( R\right) ,...,\overset{d}{\underset{j=1}{\Sigma }}%
u_{j}^{k}\left( R\right) \right) \right] ^{\frac{1}{\min \left\{
p_{1},...,p_{d}\right\} -1}}\overset{d}{\underset{i=1}{\Sigma }}A_{i}\left(
R\right) }{\overset{d}{\underset{i=1}{\Sigma }}u_{i}^{k}\left( R\right) }%
,k\geq 1.  \label{28}
\end{equation}%
Moreover, taking into account the monotonicity of 
\begin{equation*}
\left( \overset{d}{\underset{i=1}{\Sigma }}u_{i}^{k}\left( R\right) \right)
_{k\geq 1},
\end{equation*}%
there exists 
\begin{equation*}
L\left( R\right) :=\lim_{k\rightarrow \infty }\overset{d}{\underset{i=1}{%
\Sigma }}u_{i}^{k}\left( R\right) .
\end{equation*}%
We prove that $L\left( R\right) $ is finite. Indeed, if not, we let $%
k\rightarrow \infty $, in (\ref{28}) and the assumption (\ref{10}) leads us
to a contradiction. Since $u_{i}^{k}\left( R\right) $ are increasing
functions, it follows that the map $L:\left( 0,\infty \right) \rightarrow
\left( 0,\infty \right) $ is nondecreasing and 
\begin{equation*}
\overset{d}{\underset{i=1}{\Sigma }}u_{i}^{k}\left( r\right) \leq \overset{d}%
{\underset{i=1}{\Sigma }}u_{i}^{k}\left( R\right) \leq L\left( R\right)
,\forall r\in \left[ 0,R\right] ,\forall k\geq 1.
\end{equation*}%
Thus the sequences $\left( u_{i}^{k}\left( R\right) \right)
_{i=1,...,d}^{k\geq 1}$ are bounded from above on bounded sets. We now
define the following quantities 
\begin{equation*}
u_{i}\left( r\right) :=\lim_{k\rightarrow \infty }u_{i}^{k}\left( r\right) 
\text{ for all }r\geq 0\text{ and }i=1,...,d.
\end{equation*}%
Then $u_{i}$ is a positive solution of (\ref{222}).

Next, we show that $u_{i}$ ($i=1,...,d$), is a large solution of (\ref{222}%
). Let us remark that by (\ref{233}) we have the following estimate 
\begin{equation*}
u_{i}\left( r\right) \geq \beta _{i}+\left[ f_{i}\left( \beta _{1},...,\beta
_{d}\right) \right] ^{1/\left( p_{i}-1\right) }A_{i}\left( r\right) \text{,
for all }r\geq 0\text{ and }i=1,...,d\text{.}
\end{equation*}%
It follows from the assumption $f_{i}$ are positive functions and $%
A_{i}\left( \infty \right) =\infty $, that $u_{i}$ ($i=1,...,d$) is a large
solution of (\ref{222}) and so $u_{i}$ is a positive entire large solution
of \ (\ref{11}). Thus any large solution of (\ref{222}) provides a positive
entire large solution of (\ref{11}) with $u_{i}\left( 0\right) =\beta _{i}$.
Since $\beta _{i}\in \left( 0,\infty \right) $ ($i=1,..,d$) was chosen
arbitrarily, it follows that (\ref{11}) has infinitely many positive entire
large solutions.

(ii) Assume that 
\begin{equation*}
\sup_{s\geq 0}\left[ \underset{i=1}{\overset{d}{\sum }}\left( 1+f_{i}\left(
s,...,s\right) \right) ^{\frac{1}{\min \left\{ p_{1},...,p_{d}\right\} -1}}%
\right] <\infty 
\end{equation*}%
holds, then by (\ref{28}) we have 
\begin{equation*}
L\left( R\right) :=\lim_{k\rightarrow \infty }\overset{d}{\underset{i=1}{%
\Sigma }}u_{i}^{k}\left( R\right) <\infty .
\end{equation*}%
On the other hand%
\begin{equation*}
\overset{d}{\underset{i=1}{\Sigma }}u_{i}^{k}\left( r\right) \leq \overset{d}%
{\underset{i=1}{\Sigma }}u_{i}^{k}\left( R\right) \leq L\left( R\right) ,%
\text{ }\forall r\in \left[ 0,R\right] ,\forall k\geq 1.
\end{equation*}%
So the sequences $\left( u_{i}^{k}\left( R\right) \right)
_{i=1,...,d}^{k\geq 1}$ are bounded from above on bounded sets.

Let 
\begin{equation*}
u_{i}\left( r\right) :=\lim_{k\rightarrow \infty }u_{i}^{k}\left( r\right) 
\text{ for all }r\geq 0\text{ and }i=1,...,d.
\end{equation*}%
Then $u_{i}\left( r\right) $, $\left( i=1,...,d\right) $ is a positive
solution of (\ref{222}). It follows from (\ref{26}) that $u_{i}\left(
r\right) $, $\left( i=1,...,d\right) $ is bounded, which implies that (\ref%
{11}) has infinitely many positive entire bounded solutions. This concludes
the proof of Theorem \textbf{\ref{3}}.

\begin{equation*}
\end{equation*}

\begin{remark}
If (C1), (C2), (C3) are satisfied then%
\begin{equation*}
\int_{a}^{\infty }\frac{ds}{f_{j}^{1/\left( \min \left\{
p_{1},...,p_{d}\right\} -1\right) }\left( s,...,s\right) }=\infty \text{ for
all }j=\overline{1,d}.
\end{equation*}
\end{remark}

\begin{remark}
(see \cite{CD}) If (C1)-(C2) and%
\begin{equation*}
\int_{a}^{\infty }\frac{ds}{f_{j}^{1/\left( \min \left\{
p_{1},...,p_{d}\right\} -1\right) }\left( s,...,s\right) }=\infty \text{ for
all }j=\overline{1,d},
\end{equation*}%
are satisfied, then%
\begin{equation*}
\int_{a}^{\infty }\frac{dt}{\left( \int_{0}^{t}f_{j}\left( s,...,s\right)
ds\right) ^{1/\min \left\{ p_{1},...,p_{d}\right\} }}=\infty \text{ for all }%
j=\overline{1,d}.
\end{equation*}
\end{remark}

\begin{flushright}
\begin{tabular}{l}
Address \\ 
Drago\c{s}-P\u{a}tru Covei$^{1}$ \\ 
$^{1}$Constantin Brâncu\c{s}i University of Târgu-Jiu, \\ 
Calea Eroilor, No 30, Târgu-Jiu, Gorj, \\ 
România. \\ 
covdra@yahoo.com%
\end{tabular}
\end{flushright}

\end{document}